\documentclass{amsart}
\usepackage[all]{xy}
\usepackage{fancyhdr}
\usepackage{amsmath}
\usepackage{amsfonts,amssymb}
\usepackage{latexsym}
\RequirePackage[T1]{fontenc}

\numberwithin{equation}{section}

\newtheorem{dfn}{Definition}[section]
\newtheorem{prop}[dfn]{Proposition}

\newtheorem{thm}[dfn]{Theorem}
\newtheorem{lem}[dfn]{Lemma}

\newtheorem{rmq}[dfn]{Remark}

\newcommand{\al}{{\alpha}}
\newcommand{\la}{{\lambda}}

\newcommand{\R}{\mathbb{R}}

\newcommand{\N}{\mathbb{N}}

\newcommand{\mO}{\mathcal{O}}

\newcommand{\mg}{\mathfrak{g}}
\newcommand{\mh}{\mathfrak{h}}
\newcommand{\mm}{\mathfrak{m}}
\newcommand{\mt}{\mathfrak{t}}

\author{Anton Alekseev and Damien Calaque}

\title[Quantization of symplectic dynamical $r$-matrices]
{Quantization of symplectic dynamical $r$-matrices and the quantum composition formula}

\begin{document}

\maketitle

\begin{abstract}
In this paper we quantize symplectic dynamical $r$-matrices over a possibly nonabelian base. 
The proof is based on the fact that the existence of a star-product with a nice property (called 
strong invariance) is sufficient for the existence of a quantization. We also classify such quantizations and 
prove a quantum analogue of the classical composition formula for coboundary dynamical $r$-matrices. 
\end{abstract}

\section*{Introduction}

Let $\mh\subset\mg$ be an inclusion of Lie algebras and $H\subset G$ the corresponding inclusion of Lie groups. 
Let $U\subset\mh^*$ be an invariant open subset and let $Z\in(\wedge^3\mg)^\mg$. 
A {\em (coboundary) dynamical $r$-matrix} is a $\mh$-equivariant map $r:U\to\wedge^2\mg$ satisfying the 
{\em (modified) classical dynamical Yang-Baxter equation} 
$$
\frac12[r(\la),r(\la)]-\sum_ih_i\wedge\frac{\partial r}{\partial\la^i}(\la)=Z
\quad(\la\in U)\,,
$$
where $(h_i)_i$ and $(\la^i)_i$ are dual bases of $\mh$ and $\mh^*$, respectively. Then (following \cite{X,EE1}) 
\begin{equation}\label{eq:pi}
\pi_r:=\pi_{lin}+\sum_i\frac{\partial}{\partial\la^i}\wedge\overrightarrow{h_i}+\overrightarrow{r(\la)}\,,
\end{equation}
together with $Z$, defines a $H$-invariant ($\mg$-)quasi-Poisson structure on $M=U\times G$. 
Here $\pi_{lin}$ is the linear Poisson structure on $U\subset\mh^*$. \\
Such a dynamical $r$-matrix is called {\em symplectic} if the $\mg$-quasi-Poisson manifold $(M,\pi_r,Z)$ is symplectic, 
i.e.~if $\pi_r^\#:T^*M\to TM$ is invertible and $Z=0$. \\

By a {\em dynamical twist quantization} of $r(\la)$ we mean a $\mh$-equivariant map 
$J=1\otimes1+O(\hbar):U\to\otimes^2U\mg[[\hbar]]$ satisfying the {\em semi-classical limit condition} 
$$
J(\la)-J^{2,1}(\la)=\hbar r(\la)+O(\hbar^2)\quad(\la\in U)
$$
and the {\em modified dynamical twist equation} 
$$
J^{12,3}(\la)*_{PBW}J^{1,2}(\la+\hbar h^{(3)})
=\Phi^{-1}J^{1,23}(\la)*_{PBW}J^{2,3}(\la)\quad(\la\in U)\,.
$$
where $\Phi\in(\otimes^3U\mg)^\mg[[\hbar]]$ is an associator quantizing $Z$, of which we know the existence from 
\cite[Proposition 3.10]{Dr}. 
Recall that $*_{PBW}$ is the Poincaré-Birkhoff-Witt star-product given on polynomial functions as the pull-back of 
the usual product in $U(\mh_\hbar)$\footnote{$\mh_\hbar=\mh[[\hbar]]$ with bracket $[,]_\hbar:=\hbar[,]_\mh$. } by 
the symmetrization map $~{\rm sym}\!:\!S(\mh)[[\hbar]]\to U(\mh_\hbar)$. We also made use of the following notations: 
$$
J^{12,3}(\la):=(\Delta\otimes{\rm id})(J(\la))\quad\textrm{and}
$$ 
$$
J^{1,2}(\la+\hbar h^{(3)}):=\sum_{k\geq0}\frac{\hbar^k}{k!}\sum_{i_1\cdots i_k}
\frac{\partial^kJ}{\partial\la^{i_1}\cdots\partial\la^{i_k}}(\la)\otimes h_{i_1}\cdots h_{i_k}\,.
$$
In this paper we prove the following generalization of \cite[Theorem 5.3]{X2} to the case of a nonabelian base: 
\begin{thm}\label{main}
Any symplectic dynamical $r$-matrix admits a dynamical twist quantization (with $\Phi=1$).  
\end{thm}

Furthermore, two dynamical twist quantizations $J_1,J_2:U\to\otimes^2U\mg[[\hbar]]$ of $r$ are said 
to be {\em gauge equivalent} if there exists a $\mh$-equivariant map $T=1+O(\hbar):U\to U\mg[[\hbar]]$ 
such that 
$$
T^{12}(\la)*_{PBW}J_1(\la)=J_2(\la)*_{PBW}T^1(\la+\hbar h^{(2)})*_{PBW}T^2(\la)\,.
$$
In this context we also prove the following generalization of \cite[Section 6]{X2} to the case of a nonabelian base, 
asserting that dynamical twist quantizations are classified by the {\em second dynamical $r$-matrix cohomology} 
(see Definition \ref{dfn:dyncoho}): 
\begin{thm}\label{thm:equiv}
Let $r$ be a symplectic dynamical $r$-matrix. Then the space of gauge equivalence classes of dynamical twist 
quantizations of $r$ (with $\Phi=1$) is an affine space modeled on $H^2_r(U,\mg)[[\hbar]]$. 
\end{thm}

A class of examples of symplectic dynamical $r$-matrices is given by {\em nondegenerate reductive splittings}. A reductive 
splitting $\mg=\mh\oplus\mm$ (with $[\mh,\mm]\subset\mm$) is called nondegenerate if there exists $\lambda\in\mh^*$ 
for which $\omega(\la)\in\wedge^2\mm^*$ defined by $\omega(\la)(x,y)=<\la,[x,y]_{|\mh}>$ is nondegenerate. 
Then $r_\mh^\mm(\la):=-\omega(\la)^{-1}\in\wedge^2\mm\subset\wedge^2\mg$ defines a symplectic dynamical $r$-matrix on the 
invariant open subset $\{\la|\textrm{det}\omega(\la)\neq0\}\subset\mh^*$ (see \cite[Proposition 1]{FGP} or 
\cite[Theorem 2.3]{X}; see also \cite[Proposition 1.1]{EE2} for a more algebraic proof). Moreover, one can use it 
to ``compose'' dynamical $r$-matrices (see \cite[Proposition 1]{FGP} and \cite[Proposition 0.1]{EE2}): 
\begin{prop}[The composition formula]\label{prop:classicalcomp}
Assume that $\mh=\mt\oplus\mm$ is a nondegenerate reductive splitting and let $r_\mt^\mm:\mt^*\supset V\to\wedge^2\mh$ be 
the corresponding symplectic dynamical $r$-matrix. If $\rho:\mh^*\supset U\to\wedge^2\mg$ is a dynamical $r$-matrix with 
$Z\in(\wedge^3\mg)^\mg$, then 
$$
\theta_\rho:=r_\mt^\mm+\rho_{|\mt^*}:\mt^*\supset U\cap V\longrightarrow\wedge^2\mg
$$
is a dynamical $r$-matrix with the same $Z$. 
\end{prop}
We prove that one can ``quantize'' the map $\rho\mapsto\theta_\rho$:
\begin{thm}\label{thm:quantumcomp}
With the hypothesis of Proposition \ref{prop:classicalcomp}, there exists a map 
$$
\Theta:\{\textrm{Dynamical twist quantizations of }\rho\}\longrightarrow
\{\textrm{Dynamical twist quantizations of }\theta_\rho\}
$$
which keeps the associator $\Phi$ fixed. 
\end{thm}
The paper is organized as follows. \\
In Section 1 we first recall basic facts about quasi-Poisson manifolds, compatible quantizations and 
related results. We then give a sufficient condition for the existence of dynamical twist quantizations. \\
In Section 2 we give a very short proof of Theorem \ref{main}, using a result stating the existence of a quantum momentum 
map (see also \cite{F1,Kr}) which is based on Fedosov's well-known globalization procedure \cite{F,F2}. We start the 
section with a summary of basic ingredients of Fedosov's construction. \\
In Section 3 we prove Theorem \ref{thm:equiv} using a variant of the well-established classification of star-products 
on a symplectic manifold by formal series with coefficients in the second De Rham cohomology group of the manifold. \\
In Section 4 we prove a quantum analogue of the composition formula for classical dynamical $r$-matrices. We start 
with a new proof of the classical composition formula (Proposition \ref{prop:classicalcomp}) using (quasi-)Poisson 
reduction. We then derive its quantum counterpart (Theorem \ref{thm:quantumcomp}) using quantum reduction. \\

\noindent{\bf Notations. }
We denote by $\mO_{\mh^*}=S(\mh)$ polynomial functions on $\mh^*$, $\mO_U=C^\infty(U)\supset\mO_{\mh^*}$ smooth functions 
on $U\subset\mh^*$, and $\mO_G$ the ring of smooth functions on the (eventually formal) Lie group $G$. For any element 
$x\in\mg$ we denote by $\overrightarrow{x}$ (resp. $\overleftarrow{x}$) the corresponding left 
(resp. right) invariant vector field on $G$. \\

\noindent{\bf Acknowledgements. }
We thank Benjamin Enriquez for pointing to us the problem of the quantum composition formula and for a very stimulating 
series of discussions during his short visit in Geneva. 
A.A. acknowledges the support of the Swiss National Science Foundation. 
D.C. thanks the {\em Section de Mathématiques de l'Université de Genève}, 
the {\em Institut des Hautes \'Etudes Scientifiques} and 
the {\em Max-Planck-Institut für Mathematik Bonn} where parts of this work were done. 

\section{A sufficient condition for the existence of a dynamical twist quantization}

Let $r:\mh^*\supset U\to\wedge^2\mg$ be a coboundary dynamical $r$-matrix. Denote by $\pi_r$ the bivector field on 
$M=U\times G$ given by (\ref{eq:pi}). 

\subsection{Quasi-Poisson manifolds and their quantizations}

Recall from \cite{AK,AKM} that a {\em ($\mg$-) quasi-Poisson manifold} is a manifold $X$ together with a $\mg$-action 
$\rho:\mg\to\mathfrak X(X)$, an invariant bivector field $\pi\in\Gamma(X,\wedge^2TX)^\mg$ and an element 
$Z\in(\wedge^3\mg)^\mg$ such that 
\begin{equation}\label{eq:quasi-Poisson}
[\pi,\pi]=\rho(Z)\,.
\end{equation}
Let $\{f,g\}:=<\pi,df\wedge dg>$ ($f,g\in\mO_X$) be the corresponding {\em quasi-Poisson bracket}. Then equation 
(\ref{eq:quasi-Poisson}) is equivalent to 
$$
\{\{f,g\},h\}+\{\{h,f\},g\}+\{\{g,h\},f\}=<\rho(Z),df\wedge dg\wedge dh>\quad(f,g,h\in\mO_X)\,.
$$

One can define the {\em quasi-Poisson cochain complex} of $(X,\pi,Z)$ as follows: 
$k$-cochains are $C^k_\pi(X)=\Gamma(X,\wedge^kTX)^\mg$ and the differential is $d_\pi=[\pi,-]$. The fact 
that $d_\pi\circ d_\pi=0$ follows from an easy calculation: 
$$
d_\pi\circ d_\pi(x)=[\pi,[\pi,x]]=\frac12[[\pi,\pi],x]=\frac12[\rho(Z),x]=0\quad(x\in C^k_\pi(X))\,.
$$

Let us now fix an associator $\Phi\in(\otimes^3U\mg)^\mg[[\hbar]]$ quantizing $Z$ (we know it exists from 
\cite[Proposition 3.10]{Dr}). Following \cite[Definition 4.4]{EE1}, by a {\em quantization} of a given quasi-Poisson 
manifold $(X,\pi,Z)$ we mean a series $*\in{\rm Bidiff}(X)^\mg[[\hbar]]$ of invariant bidifferential operators such that 
\begin{itemize}
\item $f*g=fg+O(\hbar)$ for any $f,g\in\mO_X$, 
\item $f*g-g*f=\hbar\{f,g\}+O(\hbar^2)$ for any $f,g\in\mO_X$, and 
\item if we write $m_*(f\otimes g):=f*g$ for $f,g\in\mO_X$, and 
$\widetilde\Phi:=S^{\otimes3}(\Phi^{-1})$, then\footnote{We thank Pavel Etingof for pointing to us 
that one has to use $\widetilde\Phi$ instead of $\Phi$ in this definition. } 
\begin{equation}
m_*\circ(m_*\otimes\textrm{id})=m_*\circ(\textrm{id}\otimes m_*)\circ\rho^{\otimes3}(\widetilde\Phi)\,.
\label{eq:quasi-quant}
\end{equation}
\end{itemize}
Here, $S$ denotes the antipode of $U\mg$. \\

One has a natural notion of {\em gauge transformation} for quantizations. It is given by an element 
$Q={\rm id}+O(\hbar)\in{\rm Diff}(X)^\mg[[\hbar]]$ that act on $*$'s in the usual way:
$$ 
f*^{(Q)}g:=Q^{-1}(Q(f)*Q(g))\quad(f,g\in\mO_X)\,.
$$
More precisely, if $(\Phi,*)$ is a quantization of $(Z,\pi)$ then $(\Phi,*^{(Q)})$ is also. In this case we say that 
$*$ and $*^{(Q)}$ are {\em gauge equivalent}. 

\subsection{Classical and quantum momentum maps}

Let $(X,\pi,Z)$ be $\mg$-quasi-Poisson manifold and let $\mathcal G$ be a Lie algebra with Lie group ${\bf G}$. 
A {\em momentum map} is a smooth $\mg$-invariant map $\mu:M\to\mathcal G^*$ such that $\mu_*\pi=\pi_{lin}$ and 
for which the corresponding infinitesimal action $\mathcal G\to\mathfrak X(X);x\mapsto\{\mu^*x,-\}$ integrates to 
a right action of ${\bf G}$. 

Let us describe the reduction procedure with respect to a given momentum map $\mu$. 
First of all ${\bf G}$ acts on $\mu^{-1}(0)$ and hence one can define the reduced space $X_{red}:=\mu^{-1}(0)/{\bf G}$. 
Let us assume that it is smooth (this is the case when $0$ is a regular value and ${\bf G}$ acts freely); its function 
algebra is $\mO_{X_{red}}=\mO_X^{\mathcal G}/(\mO_X^{\mathcal G}\cap\mathcal I_0)$, where $\mathcal I_0$ is the ideal 
generated by $\textrm{im}(\mu^*)$. Since $\mu$ is $\mg$-invariant then $\mg$ acts on $\mu^{-1}(0)$. Moreover the 
$\mg$-action and the $\mathcal G$-action commute (because $\pi$ and $\mu$ are $\mg$-invariant), consequently $\mg$ 
also acts on $X_{red}$. Now observe that $\mO_X^{\mathcal G}=\{f\in\mO_X|\{f,\mathcal I_0\}\subset\mathcal I_0\}$, 
therefore the quasi-Poisson bracket $\{,\}$ naturally induces a quasi-Poisson bracket (with the same $Z$) on 
$\mO_{X_{red}}$. In other words, $X_{red}$ inherits a structure of a quasi-Poisson manifold from the one of $X$. \\

Now assume that we are given a quantization $(*,\Phi)$ of the quasi-Poisson manifold $(X,\pi,Z)$. 
By a {\em quantum momentum map} quantizing $\mu$ we mean a map of algebras 
$$
{\bf M}=\mu^*+O(\hbar):(\mO_{\mathcal G^*}[[\hbar]],*_{PBW})\longrightarrow(\mO_X[[\hbar]],*)
$$
taking its values in $\mg$-invariant functions, and such that for any $f\in\mO_X$ and any $x\in\mathcal G$ one has 
$[{\bf M}(x),f]_*=\hbar\{\mu^*x,f\}$. 
\begin{rmq}\emph{
One only needs to know ${\bf M}$ on linear functions $x\in\mathcal G$. 
}\end{rmq}
Let us describe the quantum reduction with respect to a given quantum momentum map ${\bf M}$. 
First denote by $\mathcal I$ the right ideal generated by $\textrm{im}({\bf M})$ in $(\mO_X[[\hbar]],*)$ and observe that 
its normalizer is $\mO_X^\mh[[\hbar]]$. 
Therefore $\mO_X^\mh[[\hbar]]\cap\mathcal I$ is a two-sided ideal in $\mO_X^\mh[[\hbar]]$, and we can define the 
reduced algebra $\mathcal A:=\mO^\mh_X[[\hbar]]/\mO^\mh_X[[\hbar]]\cap\mathcal I$. 
Since $\mathcal I\cong\mathcal I_0[[\hbar]]$ then $\mathcal A\cong\mO_{X_{red}}[[\hbar]]$ 
(as $\mathbb{R}[[\hbar]]$-modules). It is easy to see that the induced product on $\mO_{X_{red}}[[\hbar]]$, together with 
$\Phi$, gives a quantization of the quasi-Poisson structure on $X_{red}$. \\

A quantization $*$ of a quasi-Poisson manifold $(M,\pi,Z)$ with a momentum map $\mu:M\to\mathcal G^*$ for which $\mu^*$ 
itself defines a quantum momentum map (it will be ${\bf M}=U(\mu^*)\circ\textrm{sym}$) is called 
{\em strongly ($\mathcal G$-)invariant}. 

\subsection{Compatible quantizations}

Let us first observe that $\pi_r$ defines a quasi-Poisson structure on $M=U\times G$. Here the action is 
$\mg\ni x\mapsto\overleftarrow{x}\in\mathfrak X(M)$ (it generates left translations). 
Remark that since $Z\in(\wedge^3\mg)^\mg$ then $\overleftarrow{Z}=\overrightarrow{Z}$. Moreover, the natural 
map $M\to\mh^*,(\la,g)\mapsto\la$ is a momentum map and the corresponding right $H$-action is given by 
$(\la,g)\cdot h:=({\rm Ad}_h^*\la,gh)$ (following the notation of the previous $\S$ we have 
$\mathcal G=\mathfrak h$). Conversely, 
\begin{prop}[\cite{X}, Proposition 2.1]
A map $\rho\in C^\infty(U,\wedge^2\mg)$ is a coboundary classical dynamical $r$-matrix if and only if 
$$
\pi=\pi_{lin}+\sum\frac{\partial}{\partial\la^i}\wedge\overrightarrow{h_i}+\overrightarrow{\rho(\la)}
$$
defines a $\mg$-quasi-Poisson structure on $U\times G$. 
\end{prop}
\begin{proof}
The proof given in \cite{X} is for the case when $Z=0$, but it admits a straightforward generalization. 
\end{proof}
Following Ping Xu (\cite{X}), by a {\em compatible quantization} of $\pi_r$ we mean a 
quantization $*'$ which is such that for any $u,v\in\mO_{\mh^*}$ and any $f\in\mO_G$, 
$u*'v=u*_{PBW}v$, $f*'u=fu$ and 
\begin{equation}\label{eq:comp}
u*'f=\sum_{k\geq0}\frac{\hbar^k}{k!}\sum_{i_1,\dots,i_k}
\frac{\partial^ku}{\partial\la^{i_1}\cdots\partial\la^{i_k}}
\overrightarrow{h_{i_1}}\cdots\overrightarrow{h_{i_k}}\cdot f\,.
\end{equation}
\begin{prop}
There is a bijective correspondence between compatible quantizations of $\pi_r$ and dynamical twist quantizations of $r$. 
\end{prop}
\begin{proof}
Let $*'$ be a compatible quantization of $\pi_r$. Since $*'$ is $G$-invariant then for all $f,g\in\mO_G$ one has 
$$
(f*'g)(\la)=\overrightarrow{J(\la)}(f,g)\quad(\la\in\mh^*)
$$
with $J:U\to\otimes^2U\mg[[\hbar]]$. Moreover 
\begin{lem}\label{lem:inv}
$*'$ is strongly $\mh$-invariant\footnote{In particular $*'$ is $H$-invariant. It was not noticed in \cite{X}, where 
the definition of compatible star-products includes this $H$-invariance property. The lemma claims that it comes for 
free (like in the classical situation). }. 
\end{lem}
\begin{proof}[Proof of the lemma. ]
Let $f=gu$ ($g\in\mO_G$ and $u\in\mO_{\mh^*}$) on $U\times G$. Then for any $h\in\mh$ one has 
\begin{eqnarray*}
h*'f-f*'h & = & h*'(gu)-(gu)*'h=h*'(g*'u)-(g*'u)*'h \\
& = & (h*'g)*'u-g*'(u*'h) \qquad\qquad(\textrm{$\Phi$ acts trivialy}) \\
& = & (g*'h+\hbar(\overrightarrow{h}\cdot g))*'u-g*'(u*'h) \\
& = & g*'([h,u]_{*'})+\hbar(\chi_{h}\cdot g)*'u=g([h,u]_{*_{PBW}})+\hbar(\chi_{h}\cdot g)u \\
& = & \hbar(g(\chi_{h}\cdot u)+(\chi_{h}\cdot g)u)=\hbar(\chi_{h}\cdot f)
\end{eqnarray*}
Hence for any $f\in\mO_M$, $h*f-f*h=\hbar(\chi_{h}\cdot f)$. 
\end{proof}
Therefore using \cite[Proposition 3.2]{X} one obtains that $J$ is $H$-equivariant. 
The following lemma ends the first part of the proof: 
\begin{lem}
$J$ satisfies the dynamical twist equation. 
\end{lem}
\begin{proof}[Proof of the lemma. ]
Let us define ${\bf L}:\mg\ni x\mapsto \overrightarrow{x}$ and ${\bf R}:\mg\ni x\mapsto \overleftarrow{x}$, 
and denote by $m^{(n)}:\mO_M^{\otimes n}\to\mO_M\,;\,f_1\otimes\cdots\otimes f_n\mapsto f_1\cdots f_n$ 
the standard $n$-fold product of functions. 
A computation in \cite{X} emphases the fact that for all $f,g,h\in\mO_G$, 
one has\footnote{The reader must pay attention to the following important remark: for any 
$P\in\otimes^n U\mg$ we denote by $\overrightarrow{P}$ (\textrm{resp.} $\overleftarrow{P}$) 
the corresponding left (\textrm{resp.} right) invariant multidifferential operator, while 
${\bf L}^{\otimes n}(P)$ (\textrm{resp.} ${\bf R}^{\otimes n}(P)$) is an element in 
$\otimes^n\textrm{Diff}(G)^{G_{left}}$ (\textrm{resp.} $\otimes^n\textrm{Diff}(G)^{G_{right}}$). 
Namely, $\overrightarrow{P}=m^{(n)}\circ\big({\bf L}^{\otimes n}(P)\big)$. }
$$
m_{*'}\circ(m_{*'}\otimes\textrm{id})(f\otimes g\otimes h)
=\overrightarrow{J^{12,3}(\la)*_{PBW}J^{1,2}(\la+\hbar h^{(3)})}(f\otimes g\otimes h)
$$
and
$$
m_{*'}\circ(\textrm{id}\otimes m_{*'})(f\otimes g\otimes h)
=\overrightarrow{J^{1,23}(\la)*_{PBW}J^{2,3}(\la)}(f\otimes g\otimes h)\,.
$$
Therefore, 
\begin{eqnarray*}
m_{*'}\circ(\textrm{id}\otimes m_{*'})\circ{\bf R}^{\otimes3}(\widetilde\Phi)(f\otimes g\otimes h)
& = & m^{(3)}\big({\bf L}^{\otimes3}(J^{1,23}(\la)*_{PBW}J^{2,3}(\la)){\bf R}^{\otimes3}(\widetilde\Phi)
(f\otimes g\otimes h\big)) \\
& = & m^{(3)}\big({\bf R}^{\otimes3}(\widetilde\Phi){\bf L}^{\otimes3}(J^{1,23}(\la)*_{PBW}J^{2,3}(\la))
(f\otimes g\otimes h)\big) \\
& = & \overleftarrow{S^{\otimes3}(\Phi^{-1})}\big({\bf L}^{\otimes3}(J^{1,23}(\la)*_{PBW}J^{2,3}(\la))
(f\otimes g\otimes h)\big) \\
& = & \overrightarrow{\Phi^{-1}}\big({\bf L}^{\otimes3}(J^{1,23}(\la)*_{PBW}J^{2,3}(\la))
(f\otimes g\otimes h)\big) \\
& = & \overrightarrow{\Phi^{-1}J^{1,23}(\la)*_{PBW}J^{2,3}(\la)}(f\otimes g\otimes h)\,,
\end{eqnarray*}
where the equality before the last one follows from the invariance of $\Phi$.  
This ends the proof of the lemma. 
\end{proof}
Conversely, let $J=\sum_\al f_{\al}A_\al\otimes B_\al$ be a dynamical twist quantization of $r$ 
($f_\al\in\mO_U[[\hbar]]$ and $A_\al,B_\al\in U\mg$). Following \cite{X} we define a $G$-invariant product 
$*'$ on $\mO_M[[\hbar]]$ by 
$$
g_1*'g_2:=\sum_{k\geq0,\al}\frac{\hbar^k}{k!}\sum_{i_1,\dots,i_k}f_\al*_{PBW}
(\overrightarrow{A_\al}\cdot\frac{\partial^kg_1}{\partial\la^{i_1}\cdots\partial\la^{i_k}})*_{PBW}
(\overrightarrow{B_\al}\overrightarrow{h_{i_1}}\dots\overrightarrow{h_{i_k}}\cdot g_2)\,.
$$
One can check by direct computations that $\mh$-equivariance of $J$ implies strong $\mh$-invariance of $*'$, and that 
the dynamical twist equation implies equation (\ref{eq:quasi-quant}). 
\end{proof}
\begin{rmq}{\em
Since $\mO_{\mh^*}=S(\mh)$ is generated as a vector space by $h^n$, $h\in\mh$ and $n\in\N$, then one can rewrite 
condition (\ref{eq:comp}) as 
$$
h^n*'f=\sum_{k=0}^n\hbar^kC_n^k(\overrightarrow{h}^k\cdot f)h^{n-k}\,.
$$
}\end{rmq}
We saw in Lemma \ref{lem:inv} that a compatible quantization always satisfies the strongly $\mh$-invariance condition. 
In what follows we show that this condition is actually sufficient for the existence of a compatible quantization. 

\subsection{A sufficient condition for the existence of a compatible quantization}\label{sec:condition}

\begin{prop}\label{prop:main}
Assume that we are given a strongly $\mh$-invariant quantization $*$ of $\pi_r$ on $M$. 
Then there exists a gauge equivalent compatible quantization $*'$ of $\pi_r$. 
Therefore there exists a dynamical twist quantization $J$ of $r$. 
\end{prop}
\begin{proof}
First observe that $h*h'-h'*h=\hbar[h,h']_\mh=[h,h']_{\mh_\hbar}$. Therefore we have an algebra morphism 
$$
a:U(\mh_\hbar)\longrightarrow(\mO_M,*)\,.
$$
Then define the algebra morphism $Q:\mO_{\mh^*\times G}=S(\mh)\otimes\mO_G\longrightarrow\mO_M$ 
as follows: 
$$
Q(fu)=f*a({\rm sym}(u))\quad(u\in S(\mh),\,f\in\mO_G)\,,
$$
where ${\rm sym}:S(\mh)[[\hbar]]\longrightarrow U(\mh_\hbar)$ is the isomorphism sending $h^n$ to $h^n$ for any 
$h\in\mh$. Thus $Q(h^n\otimes f)=f*\underbrace{h*\cdots*h}_{n~times}$, and since $*$ can be expressed as a series 
$m_0+O(\hbar)$ of bidifferential operators on $M$ then $Q$ can be expressed as a series ${\rm id}+O(\hbar)$ of 
differential operators on $M$. Moreover it is obviously $\mg$-invariant (since $*$ is), consequently we have a new 
quantization $*'$ of $\pi_r$, gauge equivalent to $*$, defined as follows: for any $f,g\in\mO_M$, 
$$
f*'g=Q^{-1}(Q(f)*Q(g))\,.
$$
Let us now check that $*'$ satisfies all Xu's properties for compatible quantizations. 
\begin{itemize}
\item for any $u,v\in S(\mh)$, 
\begin{eqnarray*}
u*'v & = & Q^{-1}\big(a({\rm sym}(u))*a({\rm sym}(v)\big) \\
 & = & Q^{-1}\big(a({\rm sym}(u){\rm sym}(v))\big) \\
 & = & Q^{-1}\big(a({\rm sym}(u*_{PBW}v))\big)=u*_{PBW}v
\end{eqnarray*}
\item let $u\in S(\mh)$ and $f\in\mO_G$, then $f*'u=Q^{-1}\big(f*a({\rm sym}(u))\big)=fu$. 
Let us now compute $u*'f$; we can assume that $u=h^n$, $h\in\mh$, and then 
\begin{eqnarray*}
u*'f & = & Q^{-1}\big(a({\rm sym}(u))*f\big)=
Q^{-1}\big(\underbrace{h*\cdots*h}_{n~times}*f\big) \\
& = & Q^{-1}\big(\sum_{k=0}^nC_n^k\hbar^k(\overrightarrow{h}^k\cdot f)*\underbrace{h*\cdots*h}_{n-k~times}\big) \\
& = & \sum_{k=0}^n\hbar^kC_n^k(\overrightarrow{h}^k\cdot f)h^{n-k}
\end{eqnarray*}
\item since $*$ is a $H$-invariant star-product, then $Q$ is a $H$-invariant gauge equivalence. 
Therefore $*'$ is also $H$-invariant. 
\end{itemize} 
The proposition is proved. 
\end{proof}
\begin{rmq}\label{rmq}\emph{
The gauge transformation $Q$ constructed above obviously satisfies $Q(h)=h$ for any $h\in\mh$. 
}\end{rmq}

\section{Quantization of symplectic dynamical $r$-matrices}

In this section we prove Theorem \ref{main}. 
We start by recalling Fedosov's construction of star-products on a symplectic manifold (for more details we refer 
to \cite{F,F2}). 

\subsection{Fedosov's star-products}

Let $(M,\omega)$ a symplectic manifold and denote by $\pi=\omega^{-1}$ the corresponding Poisson bivector. 
Then its tangent bundle $TM$ inherits a Poisson structure $\tilde\pi$ expressed locally as 
$$
\tilde\pi=\pi^{ij}(x)\frac{\partial}{\partial y^i}\wedge\frac{\partial}{\partial y^j}\,,
$$
where $y^i$'s are coordinates in the fibers. This Poisson structure is regular and constant on the symplectic leaves 
which are the fibers $T_xM$ of the bundle. Therefore it is quantized by the series of fiberwize bidifferential 
operators $\exp{(\hbar\tilde\pi)}$. It defines an associative product $\circ$ on sections of $W=\hat S(T^*M)[[\hbar]]$ 
that naturally extends to $\Omega^*(M,W)$. The center of $(\Omega^*(M,W),\circ)$ consists of forms that are constant 
in the fibers, i.e.~lying in $\Omega^*(M)[[\hbar]]$. 

By assigning the degree $2k+l$ to sections of $\hbar^kS^m(T^*M)$ there is a natural decreasing filtration 
$$
W=W_0\supset W_1\supset\cdots\supset W_i\supset W_{i+1}\supset\cdots\supset\mO_M\,.
$$

Now fix (once and for all) a torsion free connection $\nabla$ on $M$ with Christoffel's symbols $\Gamma_{ij}^k$. 
One can assume without loss of generality that it is {\it symplectic} (see \cite[Section 2.5]{F2}), which means 
that $\omega$ is parallel w.r.t.~$\nabla$. 
Then consider 
$$
\partial:\Omega^*(M,W)\to\Omega^{*+1}(M,W)
$$
its induced covariant derivative. In Darboux local coordinates we have 
$$
\partial=d+\frac1\hbar[\Gamma,-]_\circ\,,
$$
where $\Gamma=-\frac12\Gamma_{ijk}y^iy^jdx^k$ is a local $1$-form with values in $W$ 
($\Gamma_{ijk}=\omega_{il}\Gamma_{jk}^l$). One has 
$$
\partial^2=-\frac1\hbar[R,-]_\circ
$$
where $R=\frac14R_{ijkl}y^iy^jdx^k\wedge dx^l$, and $R_{ijkl}=\omega_{im}R_{jkl}^m$ is the curvature tensor of the 
symplectic connection $\nabla$. 

Let us consider more general derivations of $(\Omega^*(M,W),\circ)$ of the form 
$$
D=\partial-\delta+\frac1\hbar[r,-]_\circ
$$
where $r\in\Omega^1(M,W)$ and $\delta=\frac1\hbar[\omega_{ij}y^idx^j,-]_\circ$. 
A simple calculation yields that 
$$
D^2=-\frac1\hbar[\Omega,-]_\circ
$$
where $\Omega=\omega+R+\delta r-\partial r-\frac1\hbar r^2\in\Omega^2(M,W)$ is called the {\em Weyl curvature} of $D$. 
In particular $D$ is flat (i.e.~$D^2=0$) if and only if $\Omega\in\Omega^2(M)[[\hbar]]$ (i.e.~is a central 
$2$-form), and in this case the Bianchi identity for $\nabla$ implies that $d\Omega=D\Omega=0$. 

In computing $D^2$ one sees that $\delta:\Omega^*(M,W_k)\to\Omega^{*+1}(M,W_{k-1})$ has square zero and that the torsion 
freeness of $\nabla$ implies $\delta\partial+\partial\delta=0$. 
Then we define a homotopy operator $\kappa:\Omega^*(M,W_k)\to\Omega^{*-1}(M,W_{k+1})$ on monomials 
$a\in\Omega^p(M,S^q(T^*M))$: if $p+q\neq0$ then 
$$
\kappa(a)=\frac1{p+q}y^i\partial_{dx^i}a
$$
and otherwise $\kappa(a)=0$. One easily check that $\kappa^2=0$ and $\delta\kappa+\kappa\delta={\rm id}-\sigma$ where 
$$
\sigma:\Omega^*(M,W)\to C^\infty(M)[[\hbar]],a\mapsto a_{|dx^i=y^i=0}
$$
is the projection onto $0$-forms constant in the fibers. 
\begin{thm}[Fedosov]\label{thm:D-Omega}
For any closed $2$-form $\Omega=\omega+O(\hbar)\in\ Z^2(M)[[\hbar]]$ there exists a unique $r\in\Omega^1(M,W_3)$ 
such that $\kappa(r)=0$ and 
$$
D=\partial-\delta+\frac1\hbar[r,-]_\circ
$$
has Weyl curvature $\Omega$ and is therefore flat. 
\end{thm}
\begin{proof}
First observe that $\Omega=\omega+R+\delta r-\partial r-\frac1\hbar r^2$ with $\kappa(r)=0$ if and only if 
\begin{equation}\label{eq:r}
r=\kappa(\Omega-\omega-R+\partial r+\frac1\hbar r^2)\,.
\end{equation}

Since $\partial$ preserves the filtration and $\kappa$ raises its degree by $1$ then 
$\kappa(\Omega-\omega-R)\in\Omega^1(M,W_3)$ and the sequence $(r_n)_{n\geq3}$ defined by the iteration formula 
\begin{equation}\label{eq:recr}
r_{n+1}=r_0+\kappa(\partial r_n+\frac1\hbar r_n^2)
\end{equation}
with $r_3=\kappa(\Omega-\omega-R)$ converges to a unique element $r\in\Omega^1(M,W_3)$ which is a solution of equation 
(\ref{eq:r}). We have proved the existence. 

Conversely, for any solution $r\in\Omega^1(M,W_3)$ of (\ref{eq:r}) define $r_k=r\textrm{ mod }W_{k+1}$. 
Then $r_3=\kappa(\Omega-\omega-R)$ and the sequence $(r_n)_{n\geq3}$ satisfies (\ref{eq:recr}). 
Unicity is proved. 
\end{proof}
Such a flat derivation $D$ is called a {\em Fedosov connection (of $\nabla$-type)}. The previous theorem claims 
that they are in bijection with series $\Omega$ of closed two forms starting with $\omega$. 
\begin{thm}[Fedosov]
If $D$ is a Fedosov connection then for any $f_0\in C^\infty(M)[[\hbar]]$ there exists a unique $D$-closed section
$f\in\Gamma(M,W)$ such that $\sigma(f)=f_0$. Hence $\sigma$ establishes an isomorphism between $Z^0_D(W)$ 
and $C^\infty(M)[[\hbar]]$. 
\end{thm}
\begin{proof}
Let $f_0\in C^\infty(M)[[\hbar]]$. One has $D(f)=0$ with $\sigma(f)=f_0$ if and only if 
\begin{equation}\label{eq:f}
f=f_0+\kappa(\partial f+\frac1\hbar[r,f]_\circ)\,.
\end{equation}
Like in the proof of Theorem \ref{thm:D-Omega} we can solve (uniquely) this equation with the help of an iteration formula: 
$f_{n+1}=f_0+\kappa(\partial f_n+\frac1\hbar[r,f_n]_\circ)$. 
\end{proof}
Then $f*g=\sigma(\sigma^{-1}(f)\circ\sigma^{-1}(g))$ defines a star-product on $C^\infty(M)[[\hbar]]$ that 
quantizes $(M,\omega)$. A star-product constructed this way is called a {\em Fedosov star-product (of $\nabla$-type)} 
and is uniquely determined, once $\nabla$ is fixed, by its {\em characteristic $2$-form} 
$$
\omega_\hbar:=\frac1\hbar(\Omega-\omega)\in Z^2(M)[[\hbar]]\,.
$$
Moreover one can easily prove the following 
\begin{lem}[\cite{BCG}]\label{lem:tech}
Let $\omega^{(i)}_\hbar=\sum_{k>0}\hbar^{k-1}\omega_k^{(i)}\in Z^2(M)[[\hbar]]$ ($i=1,2$) and denote by $*_i$ the 
Fedosov star-product with characteristic two-form $\omega^{(i)}_\hbar$. 
If $\omega_\hbar^{(1)}=\omega_\hbar^{(2)}+O(\hbar^k)$ then 
$$
*^{(1)}-*^{(2)}=\hbar^{k+1}\pi^\#(\omega_k^{(1)}-\omega_k^{(2)})+O(\hbar^{k+2})\,.
$$
\end{lem}

\subsection{Fedosov's construction in the presence of symmetries}

Let $(M,\omega)$ a symplectic manifold. Let us prove two results on the compatibility of Fedosov's construction with 
group actions and hamiltonian vector fields. 
\begin{prop}[Fedosov]\label{prop:inv}
Assume that a group $G$ acts on $(M,\omega)$ by symplectomorphisms and is equipped with a $G$-invariant torsion free 
connection. Then for any $\omega_\hbar\in Z^2(M)^G[[\hbar]]$ the corresponding Fedosov star-product is $G$-invariant. 
\end{prop}
\begin{proof}
First observe that starting from a $G$-invariant torsion free connection $\nabla$ one can assume without loss of 
generality that it is symplectic (see the proof of Proposition 5.2.2 in \cite{F2}, where all expressions become 
obviously $G$-invariant). 

Then, being a symplectomorphism of $(M,\omega)$, any element $g\in G$ acts via its differential $dg$ on $(TM,\tilde\pi)$ 
as a Poisson automorphism linear in the fibers. Then its dual map $g^*:T^*M\to T^*M$ defined by 
$<g^*\xi,x>=<\xi,dg(x)>$ extends to $W$ as an automorphism. 

Finally, we need to prove that $g^*$ preserves the Fedosov connection with Weyl curvature $\Omega=\omega+\omega_\hbar$. 
On one hand the automorphism $g^*$ commutes with $\partial$ (since $\nabla$ is assumed to be $G$-invariant) and 
so $g^*R=R$. On the other hand $g^*$ also commutes with $\delta$ and $\kappa$, thus if $r$ is a solution of equation 
(\ref{eq:r}) with $\kappa(r)=0$ then so is $g^*r$. By uniqueness $g^*r=r$. We are done. 
\end{proof}
\begin{prop}[Fedosov]\label{prop:ham}
Let $H\in\mO_M$ such that $\chi=\{H,\cdot\}$ preserves a torsion free connection on $M$. Then for any 
$\omega_\hbar\in Z^2(M)[[\hbar]]$ such that $\iota_\chi\omega_\hbar=0$ the corresponding Fedosov star-product $*$ 
satisfies $H*f-f*H=\hbar(\chi\cdot f)$ for any $f\in\mO_M$. 
\end{prop}
\begin{proof}
First observe that $L_\chi\omega_\hbar=(d\iota_\chi+\iota_\chi d)\omega=0$. Therefore, the infinitesimal version of the 
previous proof ensures us that the Fedosov connection $D$ with Weyl curvature $\Omega=\omega+\omega_\hbar$ is 
$L_X$-equivariant. Hence in local Darboux coordinates it writes $D=d+\frac1\hbar[\gamma,-]_\circ$ with $L_\chi\gamma=0$, 
and $\Omega=-d\gamma-\frac1\hbar\gamma^2$. Let us compute 
\begin{eqnarray*}
D(H-\iota_\chi\gamma) & = & dH-d\iota_\chi\gamma-\frac1\hbar[\gamma,\iota_\chi\gamma]_\circ
=\iota_\chi\Omega+\iota_\chi d\gamma+\frac1\hbar[\iota_\chi\gamma,\gamma]_\circ \\
& = & \iota_\chi(\Omega+d\gamma+\frac1\hbar\gamma^2)=0
\end{eqnarray*}
Since $\sigma(H-\iota_\chi\gamma)=\sigma(H)=H$, it means that $\sigma^{-1}(H)=H-\iota_\chi\gamma$ in local Darboux 
coordinates. Consequently, for any $f\in\mO_M[[\hbar]]$ 
\begin{eqnarray*}
H*f-f*H & = & \sigma\big([H-\iota_\chi\gamma,\sigma^{-1}(f)]_\circ\big)
=\sigma\big(-[\iota_\chi\gamma,\sigma^{-1}(f)]_\circ\big) \\
& = & \sigma\big(-\iota_\chi\hbar(D-d)\sigma^{-1}(f)\big)
=\hbar\sigma\big(\iota_\chi d\sigma^{-1}(f)\big) \\
& = & \hbar\sigma\big(L_\chi\sigma^{-1}(f)-d\iota_\chi\sigma^{-1}(f)\big)
=\hbar L_\chi(f)=\hbar(\chi\cdot f)
\end{eqnarray*}
The proposition is proved. 
\end{proof}

\subsection{Proof of Theorem \ref{main}}\label{sec:proofmain}

Let $r:U\to\wedge^2\mg$ a symplectic dynamical $r$-matrix. A basis $\mathcal B$ of vector fields on $M=U\times G$ is 
given by $\mathcal B=(\dots,\partial_{\la^i},\dots,\dots,\overrightarrow{e_i},\dots)$ where $(\la^i)_i$ is a 
base of $\mh^*$ and $(e_i)_i$ is a base of $\mg$. Then one defines a torsion free connection $\nabla$ on $M$ as
$$
\nabla_{b}X=\frac12[b,X]\qquad\big(b\in\mathcal B,\,X\in\mathfrak X(M)\big)\,.
$$
Remark that $[\chi_h,b]\in{\rm span}_\R\mathcal B$ for any $b\in\mathcal B$. Therefore it follows immediately from the 
Jacobi identity that $\nabla$ is $\mh$-invariant: for all $X,Y\in\mathfrak X(M)$ and $h\in\mh$, 
$$
[\chi_h,\nabla_XY]=\nabla_{[\chi_h,X]}Y+\nabla_X[\chi_h,Y]\,.
$$
Thus from Proposition \ref{prop:ham} the Fedosov star-product $*$ with the trivial characteristic $2$-form is strongly 
$\mh$-invariant. Moreover $\nabla$ is obviously $G$-invariant, hence Proposition \ref{prop:inv} implies that $*$ is also 
$G$-invariant. 

Finally, we apply Proposition \ref{prop:main} to construct a compatible quantization of $\pi_r$. We are done. 
$\Box$

\section{Classification}

Let $r:\mh^*\supset U\to\wedge^2\mg$ a dynamical $r$-matrix. Denote by $\pi_r$ the corresponding 
$H$-invariant $\mg$-quasi-Poisson structure (\ref{eq:pi}) on $M=U\times G$ (together with $Z\in(\wedge^3\mg)^\mg$). 

\subsection{Strongly invariant equivalences and obstructions}

By a {\em strongly invariant equivalence} between two strongly $\mh$-invariant quantizations of $\pi_r$ we mean a 
$H$-invariant equivalence $Q$ (namely, $Q={\rm id}+O(\hbar)\in{\rm Diff(M)}^{G\times H}[[\hbar]]$) satisfying $Q(h)=h$ 
for any $h\in\mh\subset\mO_M$. We will now develop an analogue of the usual obstruction theory in this context. 

Let us denote by $b$ the Hochschild coboundary operator for cochains on the (commutative) algebra $\mO_M$. We start 
with the following result which is a variant of a standard one. 
\begin{prop}\label{prop:obs}
Suppose that $*_1$ and $*_2$ are two strongly invariant quantizations of $\pi_r$: 
$$
f*_ig=\sum_{k\geq0}\hbar^kC^i_k(f,g)\qquad(i=1,2)\,.
$$
Assume that $*_1$ and $*_2$ coincide up to order $n$, i.e.~$C^1_k=C^2_k$ if $k\leq n$. Then 
\begin{itemize}
\item[(1)] there exists $B\in\Gamma(M,\wedge^2TM)^{G\times H}$ and $E\in{\rm Diff}(M)^{G\times H}$ such that 
$B(h,-)=0$ and $E(h)=0$ if $h\in\mh\subset\mO_{\mh^*}$, $[\pi_r,B]=0$, and satisfying 
$$
(C_{n+1}^1-C_{n+1}^2)(f,g)=B(f,g)+(bE)(f,g)\qquad(f,g\in\mO_M)\,;
$$
\item[(2)] there exists $P\in{\rm Diff}(M)^{G\times H}$ such that $C_1=\pi_r+bP$ and $P(h)=0$ for $h\in\mh$; 
\item[(3)] if $B=[\pi_r,X]$, $X\in\mathfrak X(M)^{G\times H}$ such that $X(h)=0$, then the strongly invariant equivalence 
$Q=1+\hbar^nX+\hbar^{n+1}(E+[X,P])$ transforms $*_1$ into another strongly invariant star-product which coincides 
with $*_2$ up to order $n+1$. 
\end{itemize}
\end{prop}
\begin{proof}
We use a similar argument as in \cite{X2,BBG,BCG}. 

(1) It is well-known that $b(C^1_{n+1}-C^2_{n+1})=0$. Hence we may write 
$$
C^1_{n+1}-C^2_{n+1}=B+b(E_0)
$$ 
where $B\in\Gamma(M,\wedge^2TM)^{G\times H}$ is the skew-symmetric part of $C^1_{n+1}-C^2_{n+1}$ and 
$E_0\in{\rm Diff}(M)$. Moreover, one knows (see e.g. \cite{BCG}) that $[\pi_r,B]=0$, and it follows directly from 
the strong $\mh$-invariance property for $*_1$ and $*_2$ that $B(h,-)=0$ if $h\in\mh$. 

Since $U\times G$ admits a $G\times H$-invariant connection and $b(E_0)$ is obviously $G\times H$-invariant, then 
according to \cite[Proposition 2.1]{BBG} we can assume that $E_0$ is $G\times H$-invariant. In particular $E_0$ is 
$G$-invariant and hence $E_0(f)$, $f\in\mO_U$, is a function on $U$ only. Thus we can define a $H$-invariant vector field 
$\vec{v}$ on $U$ as follows: $<dh,\vec{v}>=E_0(h)$ for any $h\in\mh\subset\mO_U$. Now $E:=E_0-\vec{v}$ satisfies
all the required properties and $b(E)=b(E_0)-b(\vec{v})=b(E_0)$. 

(2) It is standard that $C_1=\pi_r+b(P_0)$. By repeating a similar argument as in (1) we can prove that $P_0$ can be 
chosen so that $P_0=P\in{\rm Diff}(M)^{G\times H}$ and satisfies $P(h)=0$ for any $h\in\mh$. 

The third statement (3) follows from an easy (and standard) calculation. 
\end{proof}
This proposition means that obstructions to strongly invariant equivalences are in the second cohomology group of the 
subcomplex $C^\infty(U,\wedge^*\mg)^\mh$ in the $H$-invariant quasi-Poisson cochain complex of $(M,\pi_r,Z)$. On such 
cochains $c$ the (quasi-)Poisson coboundary operator $[\pi_r,-]$ reduces to 
$d_r(c):=h_i\wedge\frac{\partial c}{\partial\lambda^i}+[r,c]$. 
\begin{dfn}\label{dfn:dyncoho}
\emph{The cohomology $H^*_r(U,\mg)$ of this cochain complex is called the 
{\em dynamical $r$-matrix cohomology} associated to $r:U\to\wedge^2\mg$. 
}\end{dfn}

\subsection{Classification of strongly invariant star-products}

Now assume that the quasi-Poisson manifold $(M,\pi_r)$ is actually symplectic and denote by $\omega_r$ the symplectic 
form; it is $G\times H$-invariant and satisfies $\iota_{\chi_h}\omega=0$ for any $h\in\mh\subset\mO_{\mh^*}$. The 
$G\times H$-invariant isomorphism 
$$
\pi_r^\#:T^*M\tilde\longrightarrow TM
$$
extends to a $G\times H$-invariant isomorphism of cochain complexes 
$$
(\Omega^*(M),d)\tilde\longrightarrow(\Gamma(M,\wedge^*TM),[\pi_r,-])
$$
which restricts to an isomorphism 
$$
(\Omega_\mh^*(M)^{G},d)\tilde\longrightarrow(C^\infty(U,\wedge^*\mg)^\mh,d_r)\,,
$$
where $\Omega^*_\mh(M):=\{\alpha\in\Omega^*(M)^H|\iota_{\chi_h}\alpha=0\,,\,\forall h\in\mh\}$. \\

Let us fix once and for all a symplectic $G\times H$-invariant connection $\nabla$ on $M$ (we know it exists) and remember 
from the previous section that for any $\omega_\hbar\in\hbar\Omega^2(M)^{G\times H}[[\hbar]]$ such that 
$d\omega_\hbar=0$ there exists a (unique) $G\times H$-invariant Fedosov star-product $*$ with characteristic $2$-form 
$\omega_\hbar$. Moreover, if $\omega_\hbar\in\Omega^2_\mh(M)[[\hbar]]$ then Proposition \ref{prop:ham} 
implies that $*$ is strongly $\mh$-invariant. Therefore we can associate a strongly invariant quantization of $\pi_r$ 
(which is actually a Fedosov star-product) to any closed two form $\omega_\hbar\in\Omega^2_\mh(M)^G[[\hbar]]$. 

In the rest of the section, all Fedosov star-products are assumed to be of $\nabla$-type and $G$-invariant (since they 
quantize the $\mg$-quasi-Poisson structure $\pi_r$). 

\begin{thm}
Two strongly invariant Fedosov star-products are equivalent by a strongly invariant equivalence if and 
only if their characteristic $2$-forms lie in the same cohomology class in $H^{G,2}_\mh(M)[[\hbar]]$. 
\end{thm}
\begin{proof}
Let $*_0$ and $*_1$ two strongly invariant Fedosov star-products with respective characteristic $2$-form 
$\omega_\hbar^{(0)}$ 
and $\omega_\hbar^{(1)}$. 

First assume that $\omega_\hbar^{(0)}=\omega_\hbar^{(1)}+d\alpha$ for some 
$\alpha=\sum_k\hbar^k\alpha^{(k)}\in\Omega^1_\mh(M)^G[[\hbar]]$, and define 
$\omega_\hbar(t)=\omega_\hbar^{(0)}+td\alpha$. Let $D_t=\partial-\delta+\frac1\hbar[r(t),-]$ be the Fedosov differential 
with Weyl curvature $\Omega(t)=\omega+\hbar\omega_\hbar(t)$. Let $H(t)\in\Gamma(M,W)$ be the solution of the equation 
$$D_tH(t)=-\alpha+\dot{r}(t)$$
with $\sigma(H(t))=0$. Then $H(t)$ is $G\times H$-invariant since $D_t$, $\alpha$ and $r(t)$ are. 
According to \cite[Theorem 5.5.3]{F2} the solution of the Heisenberg equation 
$$
\frac{dF(t)}{dt}+[H(t),F(t)]_\circ=0\,,\quad F(0)=f
$$
establishes an isomorphism $Z^0_{D_0}(W)\to Z^0_{D_1}(W),f\mapsto F(1)$ and then the corresponding series of differential 
operators $Q:(\mO_M[[\hbar]],*_0)\to(\mO_M[[\hbar]],*_1)$ is obviously $G\times H$-invariant. \\
Remember from the proof of Proposition \ref{prop:ham} that in local Darboux coordinates the Fedosov differential writes 
$D_t=d+\frac1\hbar[\gamma(t),-]_\circ$ and $\sigma_t^{-1}(h)=h-\iota_{\chi_h}\gamma(t)$ if $h\in\mh$. Now 
remark that $\dot{\gamma}(t)=\dot{r}(t)$ and that $\iota_{\chi_h}r(t)$ is independent of $t$. Hence 
$\sigma_t^{-1}(h)$ does not depend on $t$ and thus $Q(h)=h$. 

Conversely, assume that $*_0$ and $*_1$ are related by a strongly invariant equivalence with 
$[\omega_\hbar^{(0)}]\neq[\omega_\hbar^{(1)}]$ in $H^{G,2}_\mh(M)[[\hbar]]$. 
Write $\omega_\hbar^{(i)}=\sum_{k>0}\hbar^{k-1}\omega^{(i)}_k$ ($i=0,1$) and denote by $l$ the lowest integer for which 
$[\omega^{(0)}_l]\neq[\omega^{(1)}_l]$ in $H^{G,2}_\mh(M)$. Let us then define 
$$
\omega_\hbar^{(2)}=\sum_{0<k<l}\hbar^{k-1}\omega^{(0)}_k+\sum_{k\geq l}\hbar^{k-1}\omega^{(1)}_k
$$
and $*_2$ the strongly invariant Fedosov star-product with characteristic $2$-form $\omega_\hbar^{(2)}$. Since 
$[\omega_\hbar^{(2)}]=[\omega_\hbar^{(1)}]\in H^{G,2}_\mh(M)[[\hbar]]$ then $*_2$ is equivalent to $*_1$, 
and hence to $*_0$, by a strongly invariant equivalence. Then we deduce from Lemma \ref{lem:tech} that 
$C^2_{l+1}-C^0_{l+1}=\pi_r^\#(\omega^{(1)}_l-\omega^{(0)}_l)$, where $C^i_k$'s ($i=0,2$) are 
the cochains defining $*_i$. Thus it follows from Proposition \ref{prop:obs} that $\omega^{(1)}_l-\omega^{(0)}_l$ 
is exact and we obtain a contradiction. 
\end{proof}
\begin{thm}\label{thm:classistrong}
Any strongly invariant quantization of $\pi_r$ is equivalent to a strongly invariant Fedosov star-product by a 
strongly invariant equivalence. Therefore the set of strongly invariant star-products quantizing $\pi_r$ up to strongly 
invariant equivalences is an affine space modeled on $H^{G,2}_\mh(M)[[\hbar]]=H^2_r(U,\mg)$. 
\end{thm}
\begin{proof}
We follow \cite[Proposition 4.1]{BCG}. 

Let $*$ be an arbitrary strongly invariant quantization of $\pi_r$. Denote by $*_0$ the strongly 
invariant Fedosov star-product with the trivial characteristic $2$-form, which coincides with $*$ up to order $0$. 
Moreover, the skew-symmetric part of the first order term in $*-*_0$ vanishes, hence it follows from 
Proposition \ref{prop:obs} that there exists a strongly invariant equivalence $Q^{(0)}=1+\hbar Q_0$ that 
transforms $*$ into a new strongly invariant quantization $*^{(0)}$ which coincides with $*_0$ up to order $1$. 
Now the skew-symmetric part of the second order term in $*^{(0)}-*_0$ yields a closed form 
$\omega_1\in\Omega^2_\mh(M)^G$. 

Denote by $*_1$ the strongly invariant Fedosov star-product with characteristic $2$-form $\omega_1$. 
Lemma \ref{lem:tech} tells us that the skew-symmetric part of the second order term in $*^{(0)}-*_1$ vanishes, hence it  
follows from Proposition \ref{prop:obs} that there exists a strongly invariant equivalence $Q^{(1)}=1+\hbar^2 Q_1$ that 
transforms $*^{(0)}$ into a new strongly invariant quantization $*^{(1)}$ which coincides with $*_1$ up to order $2$. 

Repeating this procedure we get a sequence of strongly invariant equivalences $Q^{(k)}=1+\hbar^{k+1}Q_{k}$ ($k\geq0$) 
and a sequence of closed forms $\omega_k$ ($k>0$) such that the strongly invariant quantization $*^{(k)}$ obtained from 
$*$ by applying successively $Q^{(0)},\dots,Q^{(k)}$ coincides up to order $k+1$ with the strongly invariant Fedosov 
star-product $*_k$ with characteristic $2$-form $\omega_1+\cdots+\hbar^{k-1}\omega_k$. 

Finally, the strongly invariant equivalence $Q:=\cdots Q^{(2)}Q^{(1)}Q^{(0)}$ transform $*$ into the strongly 
invariant Fedosov star-product with characteristic $2$-form $\omega_\hbar=\sum_{k>0}\hbar^{k-1}\omega_k$. 
\end{proof}

\subsection{Classification of dynamical twist quantizations}

Let $T$ be a gauge equivalence of dynamical twist quantizations $J_1$ and $J_2$ of $r$. 
One can view $T$ as an element in ${\rm Diff}(U\times G)^{G\times H}[[\hbar]]$ such that $T(u)=u$ for any 
$u\in\mO_{\mh^*}$. Moreover if we denote by $*_i'$ the compatible quantization of $\pi_r$ corresponding to $J_i$ 
($i=1,2$) then it follows from an easy calculation that 
$$
T(f*_1'g)=T(f)*_2'T(g)\,.
$$

Conversely, any $G\times H$-invariant gauge equivalence $T$ from $*_1'$ to $*_2'$ which is such that 
$T(u)=u$ for any $u\in\mO_{\mh^*}$, that we will call from now a {\em compatible equivalence}, obviously 
gives rise to a gauge equivalence of the dynamical twist quantizations $J_1$ and $J_2$. 

Therefore, the set of dynamical twist quantization of $r$ up to gauge equivalences is in bijection with the 
set of compatible quantizations of $\pi_r$ up to compatible equivalences. \\

Remember from Proposition \ref{prop:main} and Remark \ref{rmq} that any strongly invariant quantization is 
equivalent to a compatible one by a strongly invariant equivalence. Moreover the PBW star-product has the following 
nice property: for any $h\in\mh$, $h^{*_{PBW}n}=h^n$. Hence any strongly invariant equivalence between two 
compatible quantizations is actually a compatible equivalence. Consequently: 
\begin{prop}\label{prop:invcomp}
There is a bijection 
$$
\frac{\{\textrm{strongly invariant quantizations of }\pi_r\}}{\textrm{strongly invariant equivalences}}
\longleftrightarrow 
\frac{\{\textrm{compatible quantizations of }\pi_r\}}{\textrm{compatible equivalences}}
$$
\end{prop}
~\\
\begin{proof}[End of the proof of Theorem \ref{thm:equiv}. ]
Assume that the dynamical $r$-matrix is symplectic. Then Theorem \ref{thm:equiv} follows from Proposition 
\ref{prop:invcomp} and Theorem \ref{thm:classistrong}. 
\end{proof}

\section{The quantum composition formula}

In this section we assume that $\mh=\mt\oplus\mm$ is a nondegenerate reductive splitting and we denote by 
$r_\mt^\mm:\mt^*\supset V\to\wedge^2\mh$ the corresponding symplectic dynamical $r$-matrix. 
Let $p:\mh\to\mt$ be the $\mt$-invariant projection along $\mm$. For any function $f$ on $\mh^*$ with values in 
a $\mh$-module $L$ we write $f_{|\mt^*}$ for the function $f\circ p^*$ on $\mt^*$ with values in $L$ viewed as a 
$\mt$-module; in particular if $f$ is $\mh$-invariant then $f_{|\mt^*}$ is $\mt$-invariant. 

\subsection{The classical composition formula (proof of Proposition \ref{prop:classicalcomp})}

Let $\rho:U\to\wedge^2\mg$ be a dynamical $r$-matrix with $Z\in(\wedge^3\mg)^\mg$. Then 
$\pi:=\pi_{r_\mt^\mm}+\pi_\rho$ defines a $\mg$-quasi-Poisson structure (with the same $Z$) on the manifold 
$X=V\times H\times U\times G$ which is 
\begin{itemize}
\item $H$-invariant with respect to left multiplication on $H$, 
\item $H$-invariant with respect to the right action on $U\times G$, 
\item $T$-invariant with respect to the right action on $V\times H$. 
\end{itemize}

The right diagonal $H$-action, given by $(\tau,x,\la,y)\cdot q=(\tau,q^{-1}x,\textrm{Ad}^*_q\la,yq)$, actually comes 
from a momentum map: 
\begin{eqnarray*}
\mu:X & \longrightarrow & \mh^* \\
(\tau,x,\la,y) & \longmapsto & \la-\textrm{Ad}^*_{x^{-1}}(p^*\tau)\,.
\end{eqnarray*}

Consequently we can apply the reduction with respect to $\mu$. The right $H$-invariant smooth map 
\begin{eqnarray*}
\psi:X=V\times H\times U\times G & \longrightarrow & M:=U\cap V\times G \\
(\tau,x,\la,y) & \longmapsto & (\tau,yx)
\end{eqnarray*}
restricts to a diffeomorphism $\mu^{-1}(0)/H\to M$ with inverse given by 
$$
(\tau,y)\longmapsto\overline{(\tau,1,p^*\tau,y)}\quad(\tau\in U\cap V\,,\,y\in G)\,.
$$
\begin{rmq}\emph{
From an algebraic viewpoint, we have an injective map of commutative algebras $\psi^*:\mO_M\to\mO_X$ with values in 
$\mO_X^\mh=\mO_{X/H}$ and such that, composed with the projection 
$\mO_X^\mh\to\mO_X^\mh/(\mO_X^\mh\cap<\textrm{im}(\mu^*)>)=\mO_{\mu^{-1}(0)/H}$, it becomes an isomorphism. 
}\end{rmq}
Since $\psi$ is obviously left $G$-invariant then it remains to show that the induced $\mg$-quasi-Poisson structure 
on $M$ is $\pi_{\rho_{|\mt^*}+r_\mt^\mm}$. Let $t,t'\in\mt\subset\mO_{\mt^*}$ and $f,g\in\mO_G$. 
First of all we have
$$
\{\psi^*t,\psi^*t'\}_X=\{t,t'\}_X=[t,t']=\psi^*[t,t']\,,
$$
hence $\{t,t'\}_M=[t,t']$. Then 
$$
\{\psi^*t,\psi^*f\}_X=\{t,f(yx)\}_X
=\overrightarrow{t}^H\cdot(f(yx))=(\overrightarrow{t}\cdot f)(yx)=\psi^*(\overrightarrow{t}\cdot f)\,.
$$
The third equality follows from the left $H$-invariance of $\overrightarrow{t}^H$. Thus 
$\{t,f\}_M=\overrightarrow{t}\cdot f$. Finally 
\begin{eqnarray*}
\{\psi^*f,\psi^*g\}_X(\tau,x,\la,y) & = & 
\overrightarrow{r_\mt^\mm(\tau)}^H\cdot(f(yx),g(yx))+\overrightarrow{\rho(\la)}^G\cdot(f(yx),g(yx)) \\
& = & \big(\overrightarrow{r_\mt^\mm(\tau)}\cdot(f,g)\big)(yx)
+\big(\overrightarrow{\rho(\textrm{Ad}^*_x\la)}\cdot(f,g)\big)(yx)\,.
\end{eqnarray*}
Therefore, when restricting to $\mu^{-1}(0)$ one obtains 
\begin{eqnarray*}
\{\psi^*f,\psi^*g\}_X(\tau,x,\textrm{Ad}^*_{x^{-1}}(p^*\tau),y) & = & 
\big(\overrightarrow{r_\mt^\mm(\tau)}\cdot(f,g)\big)(yx)+\big(\overrightarrow{\rho(p^*\tau)}\cdot(f,g)\big)(yx) \\
& = & \psi^*\big(\overrightarrow{(r_\mt^\mm+\rho_{|\mt^*})}\cdot(f,g)\big)\,.
\end{eqnarray*}
Therefore $\{f,g\}_M=\overrightarrow{(r_\mt^\mm+\rho_{|\mt^*})}\cdot(f,g)$. 
This ends the proof of Proposition \ref{prop:classicalcomp}. $\Box$

\subsection{Quantization of the momentum map $\mu$}

Let us first consider $(V\times H,\pi_{r_\mt^\mm})$. There is a momentum map 
\begin{eqnarray*}
\nu:V\times H & \longrightarrow & \mh^* \\
(\tau,x) & \longmapsto & -\textrm{Ad}_{x^{-1}}(p^*\tau)
\end{eqnarray*}
with corresponding right $H$-action on $V\times H$ given by $(\tau,x)\cdot q=(\tau,q^{-1}x)$. 

Like in subsection \ref{sec:proofmain} one has a $T$-invariant and $H$-invariant torsion free connexion on $V\times H$, 
therefore from Proposition \ref{prop:ham} the corresponding Fedosov star-product $*$ is both strongly $\mh$-invariant and 
strongly $\mt$-invariant\footnote{Remind that we also have a momentum map $V\times H\to\mt^*;(\tau,x)\mapsto\tau$ with 
corresponding right $T$-action given by $(\tau,x)\cdot b=(\textrm{Ad}^*_b(\tau),xb)$. }. 
 
Then Proposition \ref{prop:main} tells us that there exists a strongly $\mt$-invariant (and $H$-invariant) 
equivalence $Q$ such that $*':=*^{(Q)}$ is a compatible quantization of $\pi_{r_\mt^\mm}$. Consequently we can define the 
following algebra morphism: 
$$
{\bf N}:=Q^{-1}\circ U(\nu^*)\circ{\rm sym}:(\mO_{\mh^*}[[\hbar]],*_{PBW})\to(\mO_{V\times H}[[\hbar]],*')\,.
$$
It is obviously a quantization of the Poisson map $\nu$ and, moreover, for any $h\in\mh$ and any $f\in\mO_{V\times H}$ 
one has 
$$
[{\bf N}(h),f]_{*'}=Q^{-1}\big([\nu^*h,Q(f)]_*\big)=Q^{-1}(\hbar\{\nu^*h,Q(f)\})=\hbar\{\nu^*h,f\}\,. 
$$
In other words, ${\bf N}$ is a quantum momentum map quantizing $\nu$. \\

Let us now assume that we know a dynamical twist quantization $J(\la):U\to\otimes^2U\mg[[\hbar]]$ of $\rho(\la)$ 
(with some associator $\Phi$) and denote by $*_J$ the corresponding compatible quantization of $\pi_\rho$ on $U\times G$. 
Together with $*'$ it induces a quantization $*_J'$ of $\pi_{r_\mt^\mm}+\pi_\rho$ on $X$ (with the same $\Phi$). 
\begin{rmq}\emph{
Actually $*_J'$ is the compatible quantization corresponding to the dynamical twist quantization 
${\bf J}(\tau,\la):=
J_\mt^\mm(\tau)J(\la):(\mt\oplus\mh)^*\supset V\times U\longrightarrow\otimes^2U(\mh\oplus\mg)[[\hbar]]$ 
of the dynamical $r$-matrix 
${\bf r}(\tau,\la):=r_\mt^\mm(\tau)+\rho(\la):V\times U\longrightarrow\wedge^2(\mh\oplus\mg)$. 
Here $J_\mt^\mm$ is the dynamical twist quantizing $r_\mt^\mm$. 
}\end{rmq}
For any $f\in\mO_{\mh^*}$ we define 
${\bf M}(f):=({\bf N}\otimes{\rm inc})\circ\Delta(f)\in(\mO_{V\times H}\otimes\mO_{U\times G})[[\hbar]]=\mO_X[[\hbar]]$. Here 
${\rm inc}:\mO_{\mh^*}\hookrightarrow\mO_{U\times G}$ is the natural inclusion and 
$\Delta:\mO_{\mh^*}\to\mO_{\mh^*}\otimes\mO_{\mh^*}=\mO_{\mh^*\times\mh^*}$ is defined by 
$\Delta(f)(\la_1,\la_2)=f(\la_1+\la_2)$. 
\begin{prop}
The algebra morphism 
$$
{\bf M}:(\mO_\mh[[\hbar]],*_{PBW})\longrightarrow(\mO_X[[\hbar]],*_J')
$$
is a quantum momentum map quantizing $\mu$. 
\end{prop}

\subsection{Quantization of the composition formula (proof of Theorem \ref{thm:quantumcomp})}

Let us assume that $J$ is a dynamical twist quantization of $\rho$ and keep the notations of the previous subsection. 

Denote by $\mathcal I$ the right ideal generated by $\textrm{im}({\bf M})$ in $(\mO_X[[\hbar]],*_J')$ and consider 
the reduced algebra $\mathcal A:=\mO^\mh_X[[\hbar]]/\mO^\mh_X[[\hbar]]\cap\mathcal I$. 
Let $\Psi=\psi^*+O(\hbar)$ be the composition of $\psi^*:\mO_M[[\hbar]]\to\mO_X^\mh[[\hbar]]$ with the projection 
$\mO_X^\mh[[\hbar]]\rightarrow\mathcal A\cong\mO_{\mu^{-1}(0)/H}[[\hbar]]$. 
It is obviously bijective and $G$-invariant (since $\psi^*$ is), therefore it defines a quantization $\tilde*$ of the 
quasi-Poisson structure $\pi_{r_\mt^\mm+\rho_{|\mt^*}}$. 
We end the proof of Theorem \ref{thm:quantumcomp} using the following proposition: 
\begin{prop}
$\tilde*$ is a compatibe quantization. 
\end{prop}
\begin{proof}
First of all for any $u,v\in\mO_{\mt^*}$ one has 
$$
(\psi^*u)*_J'(\psi^*v)=u*_J'v=u*_{PBW}v=\psi^*(u*_{PBW}v)\,.
$$
Consequently $u\tilde*v=\Psi^{-1}(\Psi(u)\cdot_{\mathcal A}\Psi(v))=u*_{PBW}v$. 

Then let $u\in\mO_{\mt^*}$ and $f\in\mO_G$. On one hand 
$$
(\psi^*f)*_J'\psi^*u=(f(yx))*_J'u=f(yx)u=\psi^*(fu)
$$
and thus $f\tilde*u=fu$. 
On the other hand for $u=t^n$ ($t\in\mt$) one has 
$$
\psi^*(t^n)*_J'(\psi^*f)=(t^n)*_J'(f(yx))=\sum_{k=0}^n\hbar^kC_n^k\big((\overrightarrow{t}^H)^k\cdot (f(yx))\big)t^{n-k}
$$
$$
=\sum_{k=0}^n\hbar^kC_n^k(\overrightarrow{t}^k\cdot f)(yx)t^{n-k}
=\psi^*\Big(\sum_{k=0}^n\hbar^kC_n^k(\overrightarrow{t}^k\cdot f)t^{n-k}\Big)\,.
$$
Therefore $t^n\tilde*f=\sum_{k=0}^n\hbar^kC_n^k(\overrightarrow{t}^k\cdot f)t^{n-k}$. The proposition is proved. 
\end{proof}

\end{document}